\documentclass{article}

\usepackage{amsfonts,amssymb}
\newcommand{\all}{\cite{M99a, M99b, M00, M01, MSCDC}}
\newcommand{\ind}{1\kern-.4em 1}
\newcommand{\dist}{{\rm dist}}
\newcommand{\cg}{{\cal G}}

\newcommand{\supar}{\sup_{a\in A^r}}
\newcommand{\crr}{{\cal R}}
\newcommand{\implies}{\Rightarrow}
\newcommand{\eps}{\varepsilon}

\newtheorem{theorem}{Theorem}
\newtheorem{itlemma}{Lemma}[section] %number by section
\newtheorem{itproposition}[itlemma]{Proposition}
\newtheorem{itcorollary}[itlemma]{Corollary}
\newtheorem{itremark}[itlemma]{Remark}
\newtheorem{itdefinition}[itlemma]{Definition}
\newtheorem{itexample}[itlemma]{Example}

\newenvironment{lemma}{\begin{itlemma}\rm}{\end{itlemma}} %no-italics
\newenvironment{remark}{\begin{itremark}\rm}{\end{itremark}} %no-italics
\newenvironment{corollary}{\begin{itcorollary}\rm}{\end{itcorollary}}
\newenvironment{proposition}{\begin{itproposition}\rm}{\end{itproposition}}
\newenvironment{definition}{\begin{itdefinition}\rm}{\end{itdefinition}}
\newenvironment{example}{\begin{itexample}\rm}{\end{itexample}}

\newcommand{\text}[1]{\hbox{\rm \ #1\ \/}}

\newcommand{\bl}[1]{\begin{lemma}\label{#1}}
\newcommand{\br}[1]{\begin{remark}\label{#1}}
\newcommand{\bt}[1]{\begin{theorem}\label{#1}}
\newcommand{\bd}[1]{\begin{definition}\label{#1}}
\newcommand{\bp}[1]{\begin{proposition}\label{#1}}
\newcommand{\bc}[1]{\begin{corollary}\label{#1}}
\newcommand{\bfact}[1]{\begin{fact}\label{#1}}
\newcommand{\bex}[1]{\begin{example}\label{#1}}
\newcommand{\bem}[1]{\begin{example}\label{#1}}  %Yes, 2 different ones...
\newcommand{\ec}{\end{corollary}}
\newcommand{\eex}{\end{example}}
\newcommand{\eem}{\end{example}}
\newcommand{\el}{\end{lemma}}
\newcommand{\er}{\end{remark}}
\newcommand{\et}{\end{theorem}}
\newcommand{\ed}{\end{definition}}
\newcommand{\ep}{\end{proposition}}
\newcommand{\epr}{\end{proof}}
\newcommand{\bpr}{\begin{proof}}

\newcommand{\halmos}{\rule{1ex}{1.4ex}}

\newcommand{\beq}{\begin{eqnarray}}
\newcommand{\eeq}{\end{eqnarray}}
\newcommand{\beqn}{\begin{eqnarray*}}
\newcommand{\eeqn}{\end{eqnarray*}}
\newcommand{\bi}{\begin{itemize}}
\newcommand{\ei}{\end{itemize}}
\newcommand{\ben}{\begin{enumerate}}
\newcommand{\een}{\end{enumerate}}

\newcommand{\twoif}[4]{
\left\{ \begin{array}{ll}#1&#2\\#3&#4\end{array}\right.
}
\newcommand{\R}{{\mathbb R}}  %ams bold
\newcommand{\N}{{\mathbb N}}  %ams bold

\newcommand{\T}{{\cal T}}

\newcommand{\smni}{{\scriptscriptstyle \ni}}

\newcommand{\ds}{\displaystyle}
\newcommand{\noi}{\noindent}

\newcommand{\scr}{\scriptscriptstyle}

\title{Bounded-From-Below
Solutions of the Hamilton-Jacobi Equation for Optimal Control
Problems with Exit Times: \ Vanishing Lagrangians, Eikonal
Equations, and Shape-From-Shading \footnote{Supported  by NSF
Grant DMS95-00798 (H\'ector Sussmann, PI). Revised June 14, 2002.}
\footnote{Part of this work was carried out during the fall of
1999, while the author was a University and Louis Bevier Graduate
Fellow in the Rutgers University Department of Mathematics.
Another part of this work was completed during the spring of 2001,
while the author was an assistant professor at Texas A \& M
University-Corpus Christi. The author thanks Rutgers University
and Texas A \& M University-Corpus Christi for their hospitality
during this period. This paper is based in part on the author's
Ph.D. Dissertation under Professor H\'ector J. Sussmann. The
author thanks Professor Sussmann
 for suggesting the problems
addressed in this work.}
}
\author{Michael Malisoff\\
Department of Mathematics\\ Louisiana State University\\ Baton Rouge, LA
70803-4918 USA\\ $\mathtt{malisoff@math.lsu.edu}$}
\date{}

\begin{document}
\maketitle
\begin{abstract}
We study the Hamilton-Jacobi equation for undiscounted  exit time
control problems with general nonnegative Lagrangians using the
dynamic programming approach. We prove theorems characterizing the
value function as the unique bounded-from-below viscosity solution
of the Hamilton-Jacobi equation that is null on the target.  The
result applies to problems with the property that all trajectories
satisfying a certain integral condition must stay in a bounded
set. We allow problems for which the Lagrangian is not uniformly
bounded below by positive constants, in which  the hypotheses of
the known uniqueness results for Hamilton-Jacobi equations are not
satisfied. We apply our theorems to eikonal equations from
geometric optics, shape-from-shading equations from image
processing, and variants of the Fuller Problem.

\medbreak

\noi{\bf Key Words and Phrases:}\ optimal control, dynamic
programming, viscosity solutions, exit time problems

\medbreak \noi{\bf AMS Subject Classification:}\ 35F20, 49L25
\medbreak\medbreak
\end{abstract}

\section{Introduction}
Viscosity solutions form the basis for much current work in
control theory and optimization (cf. \cite{BARD, BDL97, BKS00,
DAL, RS00, So99}). In a recent series of papers (cf. \cite{M99a,
M99b, M00, M01, MSCDC}), we presented results characterizing the
value function in optimal control as the unique viscosity solution
of the corresponding Hamilton-Jacobi-Bellman equation (HJBE) that
satisfies appropriate side conditions.  These results apply to
very general classes of exit time problems with unbounded dynamics
and nonnegative Lagrangians, including H.J. Sussmann's  Reflected
Brachystochrone Problem (cf. \cite{S98, S01}) and other problems
with non-Lipschitz dynamics (cf. \cite{M99b, M01}). They also
apply to the Fuller Problem and eikonal equations where the
Lagrangians are not bounded below by positive constants and may
even vanish outside the target for some values of the control (cf.
\cite{M99a, M00, M01, MSCDC}).  In this note, we extend some
results of \cite{M99a, M00} on proper viscosity solutions of the
HJBE by characterizing the exit time value function as the unique
{\em bounded-from-below} viscosity solution of the corresponding
HJBE that is null on the target.  (Recall that properness of a
function $w:\R^N\to\R$ is the condition that $w(x)\to +\infty$ as
$||x||\to\infty$, which is a more stringent requirement than
boundedness from below.) This refinement applies to a large class
of deterministic exit time problems for which the Lagrangian is
not uniformly bounded below by a positive constant and for which
an extra affordability condition (namely, $(H_6)$ below) is also
satisfied.  We apply this result to several physical problems
studied in \cite{M00, So99}, including eikonal and
shape-from-shading equations, as well as variants of the Fuller
Problem that are not tractable using the well-known results or
using our earlier results. (For example, see \cite{So99},  which
imposes the requirement, which is not needed below, that the light
intensity $I$ for shape-from-shading satisfies $I(x)\le C<1$ for
all $x$ and some constant $C$; \cite{S00}, which considers
solutions of eikonal and shape-from-shading equations on bounded
sets; \cite{LRT, RT} for uniqueness of bounded solutions of
shape-from-shading equations; and \cite{MSCDC, So99} which impose
asymptotics, given in (\ref{still}) below,  which will not in
general be satisfied for the problems we consider here.)

Value function characterizations of this kind have been studied by
many authors for a variety of stochastic and deterministic optimal
control problems and for dynamic games. The characterizations have
been applied to the convergence of numerical schemes for
approximating value functions and differential game values with
error estimates, synthesis of optimal controls, singular
perturbation problems, asymptotics  problems, $H^\infty$- control,
and much more. See  for example  \cite{BARD, FLEM} and the
hundreds of references  in these books. For surveys of numerical
analysis applications of viscosity solutions, see \cite{BFS99,
S99}, and for uniqueness characterizations for the HJBE of {\em
discounted} exit time problems, see \cite{BARD}.  For uniqueness
characterizations for general Hamilton-Jacobi equations that do
not necessarily arise as Bellman equations, see \cite{A97, CIL92,
I85}.  For an appropriate stronger solution concept for a subclass
of problems, leading to a characterization of a maximal solution
as a unique solution, see \cite{CS99}. However, these earlier
characterizations cannot in general be applied to exit time
problems whose Lagrangians are not uniformly bounded below by
positive constants. In fact, one easily finds exit time problems
for which the Lagrangian is not bounded below by a positive
constant and for which the corresponding HJBE has more than one
bounded-from-below solution that vanishes on the target. Here is
an example from \cite{M00} where this  occurs: \bex{ex0}
\label{ex00} Choose the dynamics and Lagrangian
\begin{equation}
\label{dyno} \dot x(t)=u(t)\in [-1,1],\; \;  \ell(x,a)\equiv
L(x):=(x+2)^2\; (x-2)^2 x^2(x+1)^2 (x-1)^2,\end{equation}
respectively.  Let $v_1$ and $v_2$ denote
 the value functions for the exit time problem of bringing points
 to the  targets $\T_1=\{0\}$ and $\T_2=\{0, 2, -2\}$, respectively,
  using the data (\ref{dyno})
(cf. (\ref{value}) below).  Therefore, if we let ${\cal M}$ denote
the set of all measurable functions $u:[0,\infty)\to [-1,+1]$,
then \[ \ds v_j(x)=\inf_{u\in {\cal M}}\left\{ \int_0^{t_{\star
j}} L(\phi(s))\, ds: \, t_{\star j}<\infty, \, \phi(0)=x,\\
\dot\phi=u \text{a.e.}\right\}\; \text{for} j=1,2
\]
where  $t_{\star j}=\inf\{t\ge 0: \phi(t)\in\T_j\}$ for $j=1,2$.
One can easily check that $v_1$ and $v_2$ are both viscosity
solutions of the associated HJBE
\begin{equation}
\label{lll}
||Dv(x)||=(x+2)^2\; (x-2)^2 x^2(x+1)^2
(x-1)^2
\end{equation}
on $\R\setminus\T$ with the target $\T:=\T_1$ that vanish on $\T$.
One  checks that with the target $\T:=\T_1$, the problem satisfies
all hypotheses of the well-known theorems that characterize value
functions of exit time control problems as unique  viscosity
solutions of the HJBE that are zero on $\T$ (cf. \cite{BARD, BS92,
Sor93a}) except that the positive lower bound requirement on
$\ell$ is not satisfied.\hfill\halmos \eex \br{nrr} \label{nrrr}
One of the hypotheses we will make on the exit time problems in
the rest of this paper is that the running costs of trajectories
starting outside $\T$ and running for any positive time are always
positive (cf. condition $(H_5)$ below).  This positivity
hypothesis is not satisfied in the previous example, since the
trajectory $x(t)\equiv -1\not\in\T$ gives $\int_0^t L(x(s))\,
ds\equiv 0$ for all $t$.  On the other hand, all other hypotheses
we make in $\S$\ref{dh} below {\em do} hold for Example
\ref{ex00}. Therefore, under the set of assumptions in our
setting, condition $(H_5)$ cannot be removed.\hfill\halmos \er

This note is organized as follows.  In $\S$\ref{dh}, we introduce
the notation and hypotheses in force throughout most of the
sequel, including the definitions of the exit time HJBE, relaxed
controls and viscosity solutions. In $\S$\ref{main}, we state our
main result, and we also explain how this result improves what was
already known about viscosity solutions of the HJBE. Our results
apply to exit time problems that violate the usual positivity
condition on the Lagrangian (namely,  (\ref{strong}) below) and
that are also not tractable by means of $\all$. This is followed
in $\S$\ref{ml} by statements of the main lemmas.  In
$\S$\ref{prf}, we prove our main result, and $\S$\ref{sfs} gives
physical applications, including cases that are not tractable
using the known results or any of our earlier results. This is
followed in $\S$\ref{do} by variants of our main result for
discontinuous viscosity solutions and local solutions. We conclude
in $\S$\ref{noncompact} by showing how to use the methods of
\cite{M00} to extend our results to cases where the control set is
unbounded.

\section{Definitions and  Hypotheses}
\label{dh}

This note is concerned with problems of the form \[ \text{For \
each\ } x\in\R^N, \text{\ infimize\ }
\int_0^{t_x(\beta)}\ell^r(y_x(s,\beta),\beta(s))\, ds \]
\begin{equation}
\label{problem} \text{\ over \ all\ } \beta\in {\cal A} \text{\
for\ which\ } t_x(\beta)<\infty,
\end{equation}
where $y_x(\cdot,\beta)$ is defined to be the solution of the
initial value control problem
\begin{equation}
\label{ivp}
\dot y(t)=f^r(y(t),\beta(t))\; \; \text{a.e.}, \; \; y(0)=x
\end{equation}
 for each
$x\in\R^N$ and each $\beta\in{\cal A}:=\{\! \text{measurable\
functions} [0,\infty)\to A^r\}$ for a given fixed compact metric
space $A$ and possibly unbounded nonlinear control system $f$, and
$t_x(\beta):=\inf\{t\ge 0:  y_x(t,\beta)\in \T\}$ for a given
fixed set $\T\subset\R^N$.  (Depending on $f$, some choices of $x$
could give $t_x(\beta)=+\infty$ for all $\beta$, in which case the
infimum for (\ref{problem}) is $+\infty$.) Here, $A^r$ denotes the
set of all Radon probability measures on $A$ viewed as a subset of
the dual of the set $C(A)$ of all real-valued continuous functions
on $A$, and $\cal A$ has the weak-$\star$ topology, so ${\cal A}$
is the set of relaxed controls from \cite{A83,BARD, W72}. Notice
that $\cal A$ includes all measurable $\alpha:[0,\infty)\to A$,
which can be viewed as Dirac measure valued relaxed controls, and
that $A^r$ is compact.    We also consider (\ref{problem}) for
cases where $A\subset\R^M$ is closed but not bounded, in which
case we set $A^r=A$ and \[{\cal A}:= \{{\rm measurable\ functions\
}[0,\infty)\to S^r: S\subseteq A {\rm \ compact}\}\bigcup\]
\begin{equation}
\label{noncom} \{{\rm measurable\ functions \ }[0,\infty)\to A\}
\end{equation}
which of course reduces to the usual definition of ${\cal A}$ when $A$ is
compact.
For compact $S\subseteq A$ and measurable $\alpha_n,\alpha,
m:[0,\infty)\to S^r$, we set
$
h^r(x,m):=\int_S h(x,a)dm(a)$ for $x\in\R^N$ and $h=f,\ell$
for suitable $f$ and $\ell$ specified below, and
$\alpha_n\to
\alpha$ weak-$\star$ means that for all $t\ge 0$ and
for
all
Lebesgue integrable functions $B:[0,t]\to C(S)$, we have
\begin{equation}
\label{neww}
\ds \lim_{n\to\infty}\int_0^t\int_S (B(s))(a) \, {\rm
d}(\alpha_n(s))(a)\, {\rm d}s=
\int_0^t\int_S (B(s))(a) \, {\rm d}(\alpha(s))(a)\, {\rm d}s
\end{equation}
Also, recall that $STC\T$ is the {\bf small-time controllability
condition} that \[\T\; \subseteq\;  {\rm
int}\left(\crr^\eps\right) \; \text{for\ all}\; \eps\, >\, 0,\]
where
\[\crr^\eps:=\{x\in\R^N\! :  \exists t\in[0,\eps)\; \&\; \beta
\in {\cal A} \;  \text{s.t.} \; y_x(t,\beta)\in\T\}.\]  Roughly
speaking, $STC\T$ means points near $\T$ can be brought to $\T$ in
small time. We remark for later reference that $STC\T$ is a
property of the restriction of the vector fields $f(\cdot, a)$ to
neighborhoods of $\T$.
In most of what follows, we assume the
following standing hypotheses (but see $\S$\ref{noncompact}
for
analogs for cases where the control set $A$ is not assumed to be compact):

\bi\item[]\bi
\item[$(H_1)$]
$A$ is a nonempty compact metric space.
\item[$(H_2)$]
$\T\subset\R^N$ is closed and nonempty, $STC\T$.
\item[$(H_3)$]
$f$ is continuous, and $\exists L>0$ such that
$||f(x,a)-f(y,a)||\le L||x-y||$ $\forall x,y\in \R^N$ \& $a\in A$.
\item[$(H_4)$]
$\ell:\R^N\times A\to[0,\infty)$ is continuous.
\item[$(H_5)$]
If $t\in(0,\infty)$, $\beta\in {\cal A}$, and
$x\in\R^N\setminus\T$, then $\int_0^t
\ell^r(y_x(s,\beta),\beta(s))\, ds>0$.
\item[$(H_6)$]
If $x\in\R^N$ and $\beta\in{\cal A}$ are such that
$\limsup_{s\to\infty}||y_x(s,\beta)||=\infty$, then $\int_0^\infty
\ell^r(y_x(s,\beta), \beta(s))\, ds=+\infty$. \ei\ei \br{check}
Assumptions $(H_5)$-$(H_6)$ are expressed in terms of the
trajectories, rather than the HJBE data.  From the PDE point of
view, it is desirable to be able to check all of our assumptions
directly from the data $f^r=(f_1,f_2,\ldots, f_N)$, $\ell^r$, and
$\T$ from the PDE, rather than assuming complete knowledge of the
trajectories. One set of conditions on the data implying $(H_5)$
is (i) there are constants $K>0$ and $C>0$ such that
$\ell(x,a)\ge K|x_1|^C$ for all $a\in A$ and $x=(x_1,x_2, \ldots,
x_N)\in \R^N$, and (ii) if $y\in \R^{N-1}$ and
$(0,y)\in\R^N\setminus \T$, then $0\not\in \{f_1(0,y,m): m\in
A^r\}$. Conditions (i)-(ii) ensure that there is a positive cost
assigned to staying outside $\T$ on each interval of positive
length. These conditions will hold for example in the Fuller
Problem discussed below (cf. $\S$\ref{quadd}). By using a
generalized version of ``Barb\u{a}lat's lemma'', $(H_6)$ can also
be checked from the HJBE data (cf. \cite{M02}, $\S$2).\er

Before discussing the motivation for these hypotheses, note that
by the Filippov Selection Theorem (cf. \cite{W72}), all of our
results remain true if ${\cal A}$ is replaced by $\{{\rm
measurable\  functions\ } [0,\infty)\to A\}$ throughout the
preceding definitions and hypotheses as long as the sets \[{\cal
D}(x):=\{(f(x,a),\ell(x,a)): a\in A\}\] are  convex for all $x\in
\R^N$.  This follows from the fact that if all the sets ${\cal
D}(x)$ are convex, then each relaxed control $\beta\in {\cal A}$
admits a measurable function $\alpha:[0,\infty)\to A$ for which
\[
\int_0^t h^r(y_x(s,\beta),\beta(s))\, ds=\int_0^t
h(y_x(s,\alpha),\alpha(s))\, ds\; \; \forall t\ge 0,\; \; h=f,\ell
\]

We call $\T$, $A$, $f$, and $\ell$ the {\bf target, control set,
dynamics, and Lagrangian} for the problem (\ref{problem}),
respectively.   We let $\partial S$ and $\bar S$ denote the
boundary and closure for any set $S\subseteq \R^M$, respectively.

The interpretation of our standing hypotheses is as follows.
Condition $(H_5)$ has the economic interpretation that all
movement outside the target states is costly. Notice that $(H_5)$
is less stringent than requiring  $\ell(x,a)>0$ for all $x\in\R^N$
and $a\in A$, since it could be that  points $p$ for which
$\min\{\ell(p,a): a\in A\}=0$ have the property that all inputs
immediately bring $p$ to points $x$ where $\min\{\ell(x,a): a\in
A\}>0$, which can give $(H_5)$  (cf. $\S$\ref{sfs} for problems
with this property).   The condition $(H_6)$ has the
interpretation that trajectories that go further and further from
the starting point without bound are unaffordable. In other words,
trajectories that give finite total costs over $[0,\infty)$  must
stay in some bounded set. As we show in $\S$\ref{sfs} below,
$(H_6)$ holds for a general class of shape-from-shading equations
from image processing, as well as for problems with vanishing
Lagrangians that are not tractable using the known results (cf.
$\S$\ref{quadd} below). However, $(H_6)$ does not follow from
$(H_1)$-$(H_5)$ (cf. Remark \ref{ikk} below). Finally, we recall
(cf. \cite{BARD}, Chapter 3) that $(H_3)$ guarantees that
(\ref{ivp}) admits a unique solution $y_x(\cdot,\beta)$ defined on
$[0,\infty)$ that satisfies
\begin{equation}
\label{pds}
\ds\sup_{u\in {\cal A}}||y_x(t,u)-x||\le M_x t\; \text{for\ all}\; t\in
[0,
1/M_x],
\end{equation}
where $M_x:=\sup\{||f(z,a)||: a\in A, ||z-x||\le 1\}$ if this
supremum is nonzero and $M_x=1$ otherwise.

The {\bf value function}
$v$ of (\ref{problem}) is defined by
\begin{equation}
\label{value}
v(x)=\inf\left\{\int_0^{t_x(\beta)}\ell^r(y_x(s,\beta),\beta(s))\, ds:
\beta\in{\cal A}, t_x(\beta)<\infty\right\}\; \in \; [0,\infty]
\end{equation}
(but see Remark \ref{rk1} for extensions to problems with exit costs).
This note will study viscosity solutions $w$ of the HJBE
\begin{equation}
\label{hjb} \ds\supar\left\{-f^r(x,a)\cdot
Dw(x)-\ell^r(x,a)\right\}=0,\; \; \; x\not\in \T
\end{equation}
associated with the exit problem (\ref{problem}) that satisfy the
following side condition:

\bi\item[]\bi
\item[$(SC_w)$]
$w$ is bounded-from-below, and $w\equiv 0$ on $\T$
\ei\ei
We remark that the LHS in (\ref{hjb}) equals
$\sup\{-f(x,a)\cdot Dw(x)-\ell(x,a): a\in A\}$ (cf. \cite{BARD}).
When we say that a function $w$ is
 {\bf bounded-from-below}, we mean that there is a finite constant
$b$ so that $w(x)\ge b$ for all $x$ in the domain of $w$.
In some of what follows, we use the notation
\[
H_B(x,p):=\sup_{a\in B}\{-f(x,a)\cdot p-\ell(x,a)\}\] for closed
$B\subseteq A$. {}From $(H_1)$-$(H_4)$, we know that $H_B$ is
continuous for all compact sets $B\subseteq A$. We sometimes write
$H(x,p)$ to mean $H_A(x,p)$. We also set
\[B_q(p):=\{x\in\R^N: ||x-p||< q\}\; \; \forall q>0,\;
p\in\R^N.\]
 Letting $C^1(S)$ denote the set of all real-valued
continuously differentiable functions on any open subset $S$
of a Euclidean space, the definition
of viscosity solutions can then be stated as follows:

\bd{definitionvis}
Assume  $\cg\subseteq \R^N$ is open, $S\supseteq\cg$,
and $F:\R^N\times \R^N\to\R$ and
$w: S\to\R$ are continuous.
We will say that $w$ is  a {\bf
viscosity
solution} of $F(x, Dw(x))= 0$ on $\cg$ provided  the
following
conditions hold:
\begin{itemize}
\item[$(C_1)$]
If $\gamma \in C^1(\cg)$ and $x_o\in\cg$ are such that $x_o$ is a
local minimizer of $w-\gamma$, then $F(x_o,\,  D\, \gamma(x_o))\,
\ge\, 0$.
\item[$(C_2)$]
If $\lambda\in  C^1(\cg)$ and $x_1\in\cg$ are such that $x_1$ is a
local maximizer of $w-\lambda$, then $F(x_1,\, D\, \lambda(x_1))\,
\le\, 0$.
\end{itemize}\ed
We  also use the following equivalent definition of
viscosity solutions
based on the {\bf superdifferentials}
$D^+w(x)$ and {\bf subdifferentials}  $D^-w(x)$ of $w$.  Let
$\cal G$, $S$, $F$, and $w$ be as in Definition
\ref{definitionvis}, and define
\[
\displaystyle D^+w(x)\; :=\; \left\{p\,\in\, \R^N:\; \;
\limsup_{{\cal G}\smni y\to x}
\frac{w(y)-w(x)-p\cdot(y-x)}{||x-y||}\; \le\; \,
0\right\}
\]
\[
\displaystyle D^-w(x)\; :=\; \left\{p\, \in\, \R^N : \; \;
\liminf_{{\cal G}\smni y\to x}
\frac{w(y)-w(x)-p\cdot(y-x)}{||x-y||}\;  \ge\; \,
0\right\}
\]
One checks (cf. \cite{BARD}) that conditions $(C_1)$ and $(C_2)$ are
equivalent to
\bi\item[]
\begin{itemize}
\item[$(C'_1)$]
\ \ $F(x,p)\ge 0$ for all $x\in{\cal G}$ and $p\in
D^-w(x)$
\item[$(C'_2)$]
\ \ $F(x, p)\le 0$ for all $x\in{\cal G}$ and $p\in
D^+w(x)$
\end{itemize}
\ei
 \noindent respectively.  Therefore, we  equivalently define
 viscosity solutions by saying that
 $w$ is a
viscosity solution of $F(x, Dw(x))=0$ on $\cal G$ provided
conditions  $(C'_1)$-$(C'_2)$ hold.  Our  results can also
be extended to the case of {\em discontinuous} viscosity solutions
(cf. $\S$\ref{disco} below for the definitions and extensions).

\section{Statement of Main Result and Remarks}
\label{main}
Our main result will be the following:

\bt{Theorem 1} \label{thm1} Assume $(H_1)$-$(H_6)$.  If $w:
\R^N\to \R$ is a continuous function that is a viscosity solution
of the HJBE (\ref{hjb}) on $\R^N\setminus\T$, and if $w$ satisfies
$(SC_w)$, then $w\equiv v$.\et \br{rk1} Under the standing
hypotheses $(H_1)$-$(H_6)$, if the value function $v$ is finite
and continuous on $\R^N$, then $v$ itself is a viscosity solution
of the HJBE (\ref{hjb}) on $\R^N\setminus\T$ (cf. \cite{BARD}).
Since $v$ satisfies $(SC_v)$, Theorem \ref{thm1} then
characterizes $v$ as the unique viscosity solution of the HJBE
(\ref{hjb}) on $\R^N\setminus\T$ in the class of continuous
functions $w:\R^N\to\R$ that satisfy $(SC_w)$.  The assumption
that the control set $A$ is compact can be relaxed in various ways
(cf. $\S$\ref{noncompact} below).  Also, the statement of the
theorem remains true, with minor changes in the proof, if we
replace $v$ with
\[
v_g(x)=\inf\left\{\int_0^{t_x(\beta)}\ell^r(y_x(s,\beta),\beta(s))\, ds
+g(y_x(t_x(\beta),\beta)):
\beta\in{\cal A}, t_x(\beta)<\infty\right\}
\]
for any continuous bounded-from-below final cost function $g:\R^N\to\R$,
except that the boundary condition in $(SC_w)$ that $w\equiv 0$ on $\T$
is replaced by  $w\equiv g$ on $\T$.
For extensions of Theorem \ref{thm1} to {\em discontinuous} and local
viscosity solutions with possibly unbounded control sets, see
$\S\S$\ref{do}-\ref{noncompact}.
\hfill\halmos \er

\br{rk2} Theorem \ref{thm1} applies to problems that are not
tractable by means of the standard results from \cite{BARD} or
using \cite{M99a, M99b, M00, M01, MSCDC}. For example, the
undiscounted exit time problem results of \cite{BARD, Sor93a}
require
\begin{equation}
\label{strong}
\forall \eps>0,
\; \, \exists C_\eps>0
\; \, \text{s.t.}
\; \, \ell(x,a)\ge C_\eps
\; \, \forall a\in A \; \, \& \; \, \forall x\not\in B(\T,\eps),
\end{equation}
where $\dist(x,\T):=\inf\{||x-b||: b\in \T\}$ and
$B(\T,\eps):=\{p\in\R^N: \dist(p,\T)<\eps\}$, i.e., uniform
positive lower bounds for $\ell$, outside neighborhoods of $\T$.
In particular, (\ref{strong}) does not allow $\inf_a\ell(\cdot,
a)$ to vanish at any point outside $\T$, nor does it allow control
values $a$ for which $\ell(x,a)\to 0$ as $||x||\to\infty$.
Moreover, as we saw in Example \ref{ex00} above, this condition
cannot be dropped.  The examples we consider in this paper do not
in general satisfy (\ref{strong}) (cf. $\S$\ref{sfs} below). The
results of \cite{M99a, M00} apply to exit time problems violating
(\ref{strong}) and give conditions guaranteeing that $v$ is the
unique viscosity solution of the HJBE in a certain class of
functions that are either proper (where properness of a function
$w$ means that $w(x)\to +\infty$ as $||x||\to\infty$) or that
satisfy a suitable generalized properness notion.  The results of
\cite{M99a, M00} require the positivity condition $(H_5)$, but
they do not require $(H_6)$.    In \cite{MSCDC}, uniqueness
results are given for problems that violate (\ref{strong}) but
that do satisfy
\begin{equation}
\label{still}
\int_0^\infty\ell^r(y_x(s,\beta),\beta(s))\, ds<\infty \; \;
\implies\; \;
\lim_{s\to\infty}y_x(s,\beta)\in\T\; \; .
\end{equation}
As we will show in $\S$\ref{sfs} below, Theorem \ref{thm1}
applies to physical
problems from optics and image processing and to problems violating
{\em both} (\ref{strong}) and
(\ref{still}), including variants of the Fuller Problem
(cf. \cite{M99a, M00}).
  We
remark that while the results of \cite{M99a, M00} apply to cases where
(\ref{strong}) and
(\ref{still}) both fail, the conclusions of those results
are that if the value function is proper, then it is the
unique {\em proper} solution of the HJBE satisfying appropriate
side conditions.  Since we do not need to assume properness
in Theorem \ref{thm1}, our results can be viewed as an
improvement of the results of \cite{M99a} and \cite{M00}
for cases where the extra  affordability condition $(H_6)$
is also satisfied. Notice too that $(H_6)$ can be expressed
as
\begin{equation}
\label{weaker}
\int_0^\infty\ell^r(y_x(s,\beta),\beta(s))\, ds<\infty \; \implies\;
\sup_s||y_x(s,\beta)||<\infty,
\end{equation}
which is of course less restrictive than (\ref{still}) for problems with
bounded
targets (cf. $\S$\ref{quadd} below).\hfill\halmos

\er

\section{Main Lemmas}
\label{ml}

Under our standing hypotheses $(H_1)$-$(H_6)$, one proves (cf.
\cite{BARD}) that the value function $v$ is a viscosity solution
of the HJBE (\ref{hjb}) on $\R^N\setminus\T$ when $v$ is finite
and continuous. The proof follows easily from the fact that $v$
satisfies the Dynamic Programming Principle, which asserts that
\begin{equation}
\label{dpp} v(x)\;  \; =\;  \; \inf_{\alpha\in {\cal
A}}\left\{\int_0^t \ell^r( y_x(s,\alpha), \alpha(s))\, ds\, +\,
v(y_x(t,\alpha))\right\} \; \; \; \; \forall\;
x\in\R^N\end{equation} for all $t\in[0, \inf_\alpha t_x(\alpha)[$.
Our uniqueness characterizations are based on the following
representation lemmas that say that viscosity solutions of the
HJBE (\ref{hjb}) on $\R^N\setminus\T$ satisfy analogs of
(\ref{dpp}). The proofs of these lemmas are  based on uniqueness
characterizations for finite horizon control (cf. Chapter 3 of
\cite{BARD}).

\bl{lemma0} Assume $(H_1)$-$(H_4)$ are
satisfied
and  $u\in C(\bar E)$ is a viscosity solution of $H(x,D
u(x))=0$ on $E$, where $E\subset\R^N$ is bounded and open. If we
set $\tau_q(\beta)=  \inf\{ t\ge 0 :
y_q(t,\beta) \in
\partial E\}$ for each $\beta\in{\cal A}$ and $q\in E$,
then, for all $\beta\in {\cal A}$ and $q\in E$, we have
\begin{equation}
\label{halff}
 u(q)\; \; \le\; \;
\int_0^\delta\ell^r(y_q(s,\beta),\beta(s))\, ds\; \; +\; \;
u(y_q(\delta,\beta))\end{equation} for $0 \le \delta <\tau_q(\beta)$. \el

\bl{lemma1} Assume that the standing hypotheses  $(H_1)$-$(H_4)$ hold  and
that $w\in C(\bar B)$ is a viscosity solution of the HJBE $H(x,Dw(x))=0$
on
$B$,
where $B$
is open and bounded. Set
\[
T_\delta(p)\; \; :=\; \; \inf \{\, t:
\dist(y_p(t,\alpha),\,
\partial B)\le \delta,\; \alpha\in {\cal A} \}\]
 for each $p\in B$ and $\delta>0$. Then
for any $p\in B$ and any $\delta\in ]0,
\dist(p,\partial B)/2]$, we have
\begin{equation}
\label{quadruple} \! \! \! w(p)\; \; \ge\; \; \inf_{\alpha\in{\cal
A}}\left\{\! \int_0^t\! \! \ell^r(y_p(s,\alpha),\alpha(s))\, ds\;
+\; w(y_p(t,\alpha))\! \right\}
\end{equation}
for all $t\in]0,T_\delta(p)[$. \el

\noindent  Notice for future use that we can also put
$\delta=\tau_q(\beta)$
in (\ref{halff}) when $\tau_q(\beta)<\infty$. We also need the
following consequence of the Bellman-Gronwall
Inequality and the sequential compactness of ${\cal A}$
(cf. \cite{W72}):

\bl{lemma2} Let $A$ be a compact metric space, let $\{\alpha_n\}$
be a sequence in ${\cal A}$, and let $c>0$. Assume $f:\R^N\times
A\to\R^N$ satisfies  ($H_3$). Then there exists a
subsequence of $\{\alpha_n\}$ (which we do not relabel) and an
$\alpha\in {\cal A}$ such that the following conditions hold:
\bi
\item[1.]
$\alpha_n\to\alpha$ weak-$\star$ on $[0,c]$.
\item[2.]
If $x_n\to x$ in $\R^N$, then $y_{x_n}(\cdot, \alpha_n)\to
y_x(\cdot,\alpha)$ uniformly on $[0,c]$.\ei \el

\noi Finally, we need the following variant of Barb\u{a}lat's
Lemma shown in \cite{MSCDC}. Recall (cf. \cite{MSCDC}) that a
continuous function $g:\R\to [0,\infty)$ is said to be of class
${\cal MK}$ provided that $g(0)=0$ and that $g$ is even and
strictly increasing on $[0,\infty)$. For example, $x\mapsto |x|^q$
is of class ${\cal MK}$ for all constants $q>1$.  Also, if $G$  is
any function of Sontag's Class ${\cal K}$ (cf. \cite{meagre}),
then $g(s):=G(|s|)$ is of class ${\cal MK}$. {}From \cite{MSCDC},
we recall the following: \bl{meagrelemma} Let $g$ be a function of
class ${\cal MK}$,  $\phi:[0,\infty)\to \R$ be differentiable,
$\phi'$ be Lipschitz continuous, and $\int_0^\infty g(\phi(s))\,
ds<\infty$.  Then $\ds
\lim_{s\to\infty}\phi(s)=\lim_{s\to\infty}\phi'(s)=0$. \el

\section{Proof of Main Result}
\label{prf}

The proof that $w\le v$ pointwise is a repeated application of
Lemma \ref{lemma0} that we leave to the reader (cf. \cite{M99a}
for details).  It remains to show that $w\ge v$.  We omit the
superscripts $r$ to simplify notation in some of what follows.
The proof that $w\ge v$ is similar in spirit to an argument from
\cite{M99a, M00} but with a weak-$\star$ argument and a
localization based on $(H_6)$ replacing the `strong
controllability' and properness conditions used in \cite{M99a}.
Fix $x\in\R^N\setminus\T$, a constant $\kappa>w(x)$, and an
integer $J$ for which $x\in B_J(0)$.  Set \[S_\kappa=\{x\in\R^N:
w(x)<\kappa\},\] which is open by the hypothesis that $w$ is
continuous.  Set ${\cal S}=S_\kappa\cap B_J(0)$, which is bounded
and open. For each $p\in\R^N$ and  $\beta\in{\cal A}$, set
\[
\tau_p(\beta):=\inf\{t\ge 0:  y_p(t,\beta)\in
\partial({\cal S}\setminus\T)\}.
\]
Fix \[\eps\in]0,\kappa-w(x)[.\]
Set
\[I(x,t,\alpha)\; \; :=\; \;
\int_0^t\ell(y_x(s, \alpha),\alpha(s))ds \; +\;
w(y_x(t,\alpha))\]
wherever the RHS is defined.
We also set
\[
T_\delta(p)\; =\; \inf\left\{t\ge 0: \dist\left(y_p(t,\alpha),
\partial({\cal S}\setminus\T)\right)<\delta, \alpha\in {\cal
A}\right\}
\]for all  $p\in \R^N$ and $\delta>0$, and we define
$x_1:=x$, $\tau_1:=T_1(x_1)$ when $T_1(x_1)<+\infty$, and
$\tau_1$:=10 when $T_1(x_1)=+\infty$. We can then use
(\ref{quadruple}) of Lemma \ref{lemma1} to get an $\alpha_1\in
{\cal A}$ such that \[w(x_1) \ge  I(x_1,\tau_1,\alpha_1)-\eps
/4.\] (We will always assume that $\delta$ of that lemma can be
taken to be $1$. Otherwise, replace $T_{1/k}(x_k)$ in what follows
with $T_{\delta_k}(x_k)$ for an appropriate sequence
$\delta_k\downarrow 0$.) Note that $y_{x_1}(\tau_1, \alpha_1)\in
{\cal S} \setminus \T$. By induction, we define
\begin{equation}
\label{seventeen} x_k\; \; :=\; \; y_{x_{k-1}}(\tau_{k-1},
\alpha_{k-1})\; \in\; {\cal S}\setminus\T\; \; \text{for}\; \;
k=2,3, \ldots, \; \; \; \text{where}
\end{equation}
\[
\tau_k\; :=\; \twoif{
 T_{1/k}(x_k)}{\text{if}
 T_{1/k}(x_k)<+\infty}{10^k}{\text{otherwise}}.
 \]
Since $x_k \in {\cal S} \setminus \T$, we can use (\ref{quadruple})
to get an $\alpha_k \in {\cal A}$ such that
\begin{equation}
\label{3.29} w(x_k)\; \; \ge\; \; I(x_k,\tau_k,\alpha_k)\; -\;
2^{-(k+1)}\eps\; \; \text{for\ all} \; \; k\in \N.
\end{equation}
We also set $\sigma_o=0$, $\sigma_k:=\tau_1+ \ldots +\tau_k$,
$\bar
\sigma_J=\limsup_k\sigma_k$, and, for an arbitrary ${\bar a}\in A$,
\[ {\bar \alpha}_J(s):=
\left\{ \begin{array}{ll} \alpha_1(s)&\text{if}0\le s<\sigma_1,\\
\alpha_2(s-\sigma_1)&\text{if} \sigma_1 \le s <\sigma_2,\\
\vdots&\\ \alpha_k(s-\sigma_{k-1})&\text{if}\sigma_{k-1}\le s
<\sigma_k,\\ \vdots&\\ {\bar a}&\text{if}\bar \sigma_J \le s,
\end{array}\right.
\]
with the last line used if $\bar \sigma_J < +\infty$.
(We use the subscript $J$ to indicate the choice of radius in $B_J(0)$.)
{}From
the definitions of $x_k$ and $\bar \alpha_J$, we know that
\begin{equation}
\label{3.30} y_x(s, \bar \alpha_J)\; =\; y_{x_k}(s-\sigma_{k-1},
\alpha_k)\; \in\;  {\cal S} \setminus {\cal T}\; \;
\text{when}\; \; s<{\bar \sigma}_J
\end{equation}
and
\begin{equation}
\label{twotwo}
 \int_0^{\tau_k}\ell(y_{x_k}(s,\alpha_k), \alpha_k(s))\; ds\; \;
 =\; \;
\int_{\sigma_{k-1}}^{\sigma_k}\ell(y_x(s,{\bar \alpha}_J), {\bar
\alpha}_J(s))\,  ds\ge  0  \text{for  all} k.
\end{equation}
 Reapplying
(\ref{3.29}), we therefore get
\begin{eqnarray}\label{diverge}
w(x)\! \! \! &\ge&\! \! \! \int_0^{\tau_1}\ell(y_x(s,
\bar\alpha_J),\bar\alpha_J(s))\; ds+w(x_2)-\eps/4\nonumber \\
    &\ge&\! \! \!
\int_0^{\sigma_2}\ell(y_x(s, \bar\alpha_J), \bar\alpha_J(s))\; ds+w(x_3)-
\eps\left(\frac{1}{4}\! +\! \frac{1}{8}\right)\nonumber\\
    &\ge&\! \! \! \ldots\nonumber \\
    &\ge&\! \! \! I(x,\sigma_k, \bar \alpha_J)-\frac{\eps}{2}\;
\left(1-\frac{1}{2^{k}}\right)
\; \,
\forall k\in \N.
\end{eqnarray}
By (\ref{seventeen}) and the boundedness of ${\cal S}$, we can find
$\bar x_J\in\bar{\cal S}$ and a subsequence (which we will not relabel)
for
which
$x_n\to\bar x_J$.  (We later show that $\bar x_J\in
\partial [B_J(0)]\cup\T$.)
We claim that
\begin{equation}
\label{claim}
\bar \tau_J:=\inf\{\tau_{\bar x_J}(\alpha): \alpha\in {\cal A}\}
\le \limsup_k\tau_k.
\end{equation}
To see why (\ref{claim}) holds, first let $\delta\in (0,\infty)$ be given.
Assume first
that $\bar\tau_J<\infty$.
Suppose that
for
$k$
as large as desired we had $\tau_k<\bar \tau_J-\delta$.
Passing
to a subsequence, we can assume that $\tau_k\to
z\in[0,\bar\tau_J-\delta]$.
There would then exist a sequence $\tilde\tau_k\to z$ and a control
$u\in {\cal A}$ such that
\[
\dist(y_{\bar x_J}(z,u),\partial({\cal S}\setminus\T))
\leftarrow
\dist(y_{x_k}(\tilde\tau_k,u_k),\partial({\cal S}\setminus\T))\le 1/k
\to 0,
\]
where we used the definition of the $\tau_k$'s and $u$ is a weak-$\star$
limit of the $u_k$'s
on $[0,\bar{\tau}_J-\delta]$
(cf. Lemma \ref{lemma2}).  Since $z<\bar\tau_J$, this
contradicts the
definition of $\bar \tau_J$.
If on the other hand we had $\bar\tau_J=\infty$, then we arrive at the
same contradiction by replacing $\bar\tau_J-\delta$ with an arbitrary
finite positive number in the previous argument. This establishes
the claim (\ref{claim}).

Using (\ref{claim}) and passing to
a further subsequence without relabeling, we can fix a constant
$l\in [0,+\infty]$ so that
\[
l\ge \bar \tau_J\; \; \text{and}\; \; \tau_k\uparrow l.\]
Moreover, the  estimate (\ref{pds})  for Lipschitz dynamics
easily gives $\bar\tau_J=0$
iff
$\bar x_J\in \partial ({\cal S}\setminus\T)$ (cf. \cite{M00} for details).

We now use a variant of an argument from \cite{M99a} to show that
$\bar x_J\in \partial({\cal S}\setminus\T)$.  This argument, which is
a consequence of the assumption $(H_5)$,  is
as follows.  Suppose that $\bar x_J\not\in \partial ({\cal
S}\setminus\T)$,
so $l\ge \bar \tau_J>0$.
Let $M\in (0,l)$, and let $\tilde \alpha\in
{\cal A}$ be a weak-$\star$ limit of a subsequence of the
$\alpha_{k}$'s in ${\cal A}$ on $[0,M]$, which we assume to be
the sequence itself for brevity (cf. Lemma \ref{lemma2}). We
conclude from (\ref{diverge}) that
\begin{eqnarray}
\label{kkkkaty}
 0&\leftarrow&
\int_{\sigma_{k-1}}^{\sigma_{k}\wedge \{\sigma_{k-1}+M\}}
\ell(y_x(s,\bar \alpha_J),\bar \alpha_J(s))\; ds\nonumber\\[.5em] &=&
\int_0^{\tau_{k}\wedge M} \ell(y_{x_{k}}(s,
\alpha_{k}),\alpha_{k}(s))\; ds\; \; \to\; \;  \int_0^{M}
\ell(y_{{\bar x}_J}(s, \tilde \alpha),\tilde \alpha(s))\; ds.
\end{eqnarray}
The left arrow is by the divergence test applied to the integrals
in (\ref{diverge}), since $w$ is bounded below and $\ell$ is
nonnegative.  The right arrow is justified by the argument of
\cite{M99a,M00}.

If we had $\int_0^{\bar \tau_J} \ell^r(y_{\bar x_J}(s, \tilde
\alpha),\tilde \alpha(s))\; ds>0$, then $\int_0^{G} \ell^r(y_{\bar
x_J}(s, \tilde \alpha),\tilde \alpha(s))\; ds>0$ for some
$G\in(0,\bar\tau_J)$.  Since $l\ge\bar \tau_J$,
 we would reach a contradiction by
putting $M=G$ in (\ref{kkkkaty}).  It follows that
$\int_0^{\bar \tau_J} \ell^r(y_{\bar x_J}(s, \tilde
\alpha),\tilde \alpha(s))\; ds=0$.
Since we were assuming that $\bar x_J\not\in\partial({\cal
S}\setminus\T)$, we have $\bar\tau_J>0$ and $\bar x_J\not\in\T$, so
this contradicts $(H_5)$.
Therefore, it must
have been the case that
$\bar x_J
\in
\partial ({\cal S}\setminus \T)$, as needed.
Since \begin{equation}
\label{lcirc}
\partial ({\cal S}\setminus \T)\subseteq \partial
({\cal
S}_\kappa)\cup\T
\cup \partial(B_J(0)),\end{equation} we have the following cases to
consider:

\noi {\bf Case 1:\ } If $\bar x_J\in \partial ({\cal S}_\kappa)$, then
the continuity of $w$ gives $w(\bar x_J)=\kappa$.  Using (\ref{diverge}),
the nonnegativity of $\ell$, and the fact that $\eps<\kappa-w(x)$, we
conclude
that
\[w(x)\; \ge\;   w(x_k)-\eps\; \to\;
w(\bar x_J)-\eps\; >\; \kappa-(\kappa-w(x)) \; =\;  w(x),\] which
is a contradiction.  Therefore,  $\bar x_J\not\in \partial ({\cal
S}_\kappa)$.

\noi {\bf Case 2:\ }
If $\bar x_J\in\T$, then it follows
from the controllability hypothesis
$STC\T$, the continuity of $w$, $(SC_w)$,
 and the estimate (\ref{pds})
that there exist $p\in\N$, $\tilde t>0$, and $\tilde \beta\in
{\cal A}$ that are such that
\begin{equation}
\label{tou}
w(x_p)> -\eps/4, \; \; \tilde t:=t_{x_p}(\tilde\beta)<\infty, \;
\text{and}\;
\int_0^{\tilde t} \ell(y_{x_p}(s,\tilde\beta),\tilde\beta(s))\, ds<\eps/4.
\end{equation}
Combining (\ref{diverge}) and  (\ref{tou}) now gives
\[
w(x)\ge \int_0^{t_\star}
\ell(y_x(s,\bar{\bar\alpha}),\bar{\bar\alpha}(s))\,
ds-\eps
\ge v(x)-\eps,
\]
where $\bar{\bar \alpha}$ is the concatenation of $\bar\alpha_J\lceil [0,
\sigma_{p-1}]$ followed by $\tilde \beta$, and
$t_\star:=t_x(\bar{\bar\alpha})<\infty$.
This establishes that $w(x)\ge v(x)$, by the arbitrariness of $\eps$.

\noi{\bf Case 3:\ }  Since Case 1 cannot occur, and since Case 2
gives the desired conclusion $w(x)\ge v(x)$, it follows from
(\ref{lcirc}) that we can
assume that
$\bar x_J\in \partial \left[B_J(0)\right]$ in what follows.

We may assume
$\bar\sigma_J<\infty$.  (Otherwise, in what follows, replace $\bar x_J$
with
one of the $x_k$'s for which $||x_k||\ge J-2^{-J}$
and replace $\bar\sigma_J$ with the corresponding $\sigma_{k-1}$.
This is possible since $x_k\to \bar x_J\in \partial[
B_J(0)]$.)  Notice that $w(\bar x_J)<\kappa$ and
${\bar x}_J=y_x(\bar\sigma_J, \bar\alpha_J)$.
Now repeat this procedure but with the initial value
$x$ replaced by $\bar x_J$, ${\cal S}$ replaced by
${\cal S}_\kappa\cap B_{J+1}(0)$, and $\eps$
replaced by any positive
number $\eps_1<\eps/2 \wedge [\kappa-w(\bar x_J)]$
to get
a trajectory for an input $\bar{\alpha}_{J+1}$ starting
at $\bar x_J$ which wlog reaches $\partial (B_{J+1}(0))$
at time $\bar\sigma_{J+1}<\infty$.  If we now concatenate
this result with $y_x(\cdot,\bar\alpha_J)\lceil[0,\bar\sigma_J]$,
 then we get a trajectory
that coincides with $y_x(\cdot,\bar\alpha_J)$ on $[0,
\bar\sigma_J]$ and  reaches $\partial[B_{J+1}(0)]$ in finite time
$\bar{\sigma}_J+\bar{\sigma}_{J+1}$.

This process can be repeated, with $\eps$ replaced by any positive
number $\eps_q<\eps/2^q \wedge [\kappa-w(\bar x_{J+q-1})]$ and the
starting point $x$ replaced by $\bar x_{J+q-1}$ in the $q$th
iteration of this process. We can assume
$\bar{\sigma}_{J+q}<\infty$ and that all the points $\bar
x_{J+q}=y_{\bar
{x}_{J+q-1}}(\bar{\sigma}_{J+q},\bar{\alpha}_{J+q})$ obtained lie
in $\partial \left[B_{J+q}(0)\right]$ for all $q$, by the
preceding argument.  Set
\[
\ds \bar{\bar\sigma}_q=
\bar\sigma_J+\bar\sigma_{J+1}+\ldots+\bar\sigma_q \; \; \text{and}\; \;
\bar s=\limsup_q \bar{\bar\sigma}_q
\]
Fix $\bar b\in A$.
We can then set
\[ {\hat \alpha}(s):=
\left\{ \begin{array}{ll} \bar\alpha_J(s)&\text{if}0\le
s<\bar{\bar\sigma}_J,\\
\bar{\alpha}_{J+1}(s-\bar{\bar\sigma}_J)&\text{if} \bar{\bar\sigma}_J \le
s <\bar{\bar\sigma}_{J+1},\\
\vdots&\\
\bar\alpha_{J+q}(s-\bar{\bar\sigma}_{J+q-1})&\text{if}
\bar{\bar\sigma}_{J+q-1}\le
s <\bar{\bar\sigma}_{J+q},\\ \vdots&\\
\bar b& \text{if} \bar s\le s
\end{array}\right.\]
to define an input $\hat\alpha\in {\cal A}$.
A passage to the limit as $k\to\infty$ in
(\ref{diverge}) and a summation then gives
\begin{equation}
\label{wigg}
w(x)\ge
\int_0^{\bar{\bar\sigma}_q}\ell(y_x(s,\hat\alpha),\hat\alpha(s))\,
ds
+w(\bar x_q)-2\eps \; \; \text{for}\; \; \N\ni q\ge J\; .
\end{equation}
If $\bar s$
is finite, then we get
\[
\partial\left[B_{J+q+1}(0)\right]\ni
y_{\bar
x_{J+q}}(\bar\sigma_{J+q+1},\bar\alpha_{J+q+1})=y_x(\bar{\bar\sigma}_{J+q+1},\hat
\alpha)\to
y_x(\bar
s,\hat\alpha) \; \text{as}\; q\to\infty
\]
which is impossible.
Using the fact that $w$ is bounded-from-below, a passage to the limit as
$q\to\infty$ in (\ref{wigg}) therefore gives
\begin{equation}
\label{sec}
\int_0^{\infty}\ell(y_x(s,\hat\alpha),\hat\alpha(s))\,
ds \; \; \le\; \;  w(x) \; + \text{constant} \; \; <\; \; \infty
\end{equation}
Since
\[
y_x(\bar{\bar\sigma}_{J+q+1},\hat\alpha)=y_{{\bar
x}_{J+q}}(\bar\sigma_{J+q+1},
\bar\alpha_{J+q+1})\in
\partial [B_{J+q+1}(0)]
\]
for $q=1,2,\ldots$, we also have
\begin{equation}
\label{fir}
\ds
\limsup_{s\to\infty}||y_x(s,\hat\alpha)||=\infty.
\end{equation}
But (\ref{sec})-(\ref{fir}) stand in contradiction with $(H_6)$.
Consequently, it must be the case that $\bar {x}_{J+q}\in\T$ for large
enough $q$.  By the argument above, this gives the desired
inequality $w(x)\ge v(x)$ and completes the proof.

\section{Three Applications}
\label{sfs} This section shows how Theorem \ref{thm1} applies to
exit time HJBE's that are not tractable by means of the well-known
methods, including cases where the methods of $\all$ cannot be
applied.  We also show how Theorem \ref{thm1} extends results from
\cite{M00, So99} on degenerate eikonal and shape-from-shading
equations from optics and image processing.

\subsection{Vanishing Lagrangians}
\label{quadd} Theorem \ref{thm1} can be used to give uniqueness
characterizations for HJBE's that are not tractable using $\all$
or \cite{BARD, FLEM}. For example, fix $k\ge 0$, take $N=2$, and
use the exit time data
\begin{equation}
\label{datt}
\begin{array}{l}
\T=\{(k,k)\},  \; \; A=[-1,+1]\\
f(x,y,a)=(y-k\Phi(x,y),a), \; \; \ell(x,y,a)=
x^2+k(1-|a|)^2,\end{array}
\end{equation}
where $\Phi:\R^2\to [0,1]$ is any $C^1$ function that is $1$ on
$B_{k/4}((k,k))$ and $0$ on $\R^2-B_{k/2}((k,k))$.  The physical
interpretation of this data is that $\Phi$ guarantees $STC\T$ (cf.
below), and the structure of $\ell$ penalizes inputs that are not
bang-bang. This is a generalization of the Fuller Problem exit
time problem data (cf. \cite{J97, M99a, M00, MSCDC, ZB94}), which
is the case where $k=0$ in (\ref{datt}).  Recall (cf. \cite{SP96})
that the Fuller Problem admits a cost-minimizing control $\beta_z$
for each initial state $z\in \R^2$, which is defined as follows.
Set
\[
\zeta:=\{(x_1,x_2): |x_1|=Cx^2_2, \, x_1x_2\le 0\}\subset\R^2,\]
set $\zeta^{\pm}=\{(x_1,x_2)\in \zeta: \pm x_1>0\}$, and
let $A^-$ and $A^+$ denote the regions lying above and below $\zeta$
respectively, where $C>0$ is the constant root specified in \cite{SP96}.
Define the feedback $k:\R^2\to[-1,+1]$ by
$k(q)=-1$ if $q\in A^-\cup \zeta^-$,
$k(q)=1$ if $q\in A^+\cup \zeta^+$,
and $k(0,0)=0$, and let $\gamma_z$ be the closed-loop trajectory
for the feedback $k$ starting at $z$.  We then take
$\beta_z(t)=+1$ if $\gamma_z(t)\in A^+$,
$\beta_z(t)=-1$ if $\gamma_z(t)\in A^-$,
and $\beta_z(t)=0$ if $\gamma_z(t)=(0,0)$.
Let $v_k$ denote the value function
(\ref{value}) for the
exit time problem with data (\ref{datt}).

As shown in \cite{MSCDC} (see also \cite{So99}), the value
function $v=v_o$ for the Fuller Problem is the unique
bounded-from-below viscosity solution of the corresponding HJBE on
$\R^2\setminus\{0\}$ in the class of all continuous functions
$w:\R^2\to\R$ that are null at $(0,0)$. This result uses the fact
that the Fuller Problem satisfies (\ref{still}). On the other
hand, for $k> 0$, the exit time data (\ref{datt}) violate {\em
both}  (\ref{strong}) {\em and} (\ref{still}).  For example,
(\ref{strong}) is violated since $\ell(0,p,1)\equiv 0$, even
though $(0,0)\not\in\T$. Therefore, the data (\ref{datt}) is not
tractable using \cite{BARD, BS92, Sor93a}.

To see why (\ref{still}) fails for $k>0$,
let $y^k_{\scr q}(\cdot, \alpha)$ denote the
trajectory for the data (\ref{datt}),  the control $\alpha$,
 and the initial
position $q$.
For $n\in \N$ and $\beta_z$ as defined above, let
$
p(n):=(1/(2n^2),
1/n)=y^o_{\scr (0,0)}(1/n,\alpha\equiv 1)$ and
$t_n:=\inf\{t\ge 0: y^o_{\scr p(n)}(t,\beta_{p(n)})=(0,0)\}$.
Using \cite{ZB94}, we have $M:=\sup\{t_n: n\in \N\}<\infty$.  Let
$\beta$ denote the concatenation of $\beta_{p(1)}\lceil [0,t_1]$ followed
by $\alpha\equiv 1\lceil [0,1/2]$ followed by
$\beta_{p(2)}\lceil [0,t_2]$ followed by $\alpha\equiv 1\lceil[0,1/3]$
followed
by
$\beta_{p(3)}\lceil [0,t_3]$ followed by $\alpha\equiv 1\lceil[0,1/4]$
and so
on.  Since the norm of the first
coordinate of $y^o_{\scr p(n)}(\cdot, \beta_{p(n)})$
is always below $1/n^2$ (cf. \cite{SP96}),
$v_o(p(n))\le M/n^4$ for all $n$.  For $n\ge 2$, set
\[\ds \tilde t_n=\sum_{j=1}^{n-1}\left[t_j+(j+1)^{-1}\right] \; \;
\text{and}\; \;
\gamma_n(s)=\beta(s+\tilde t_n),\] so
$p(n)=y^o_{\scr p(1)}(\tilde t_n,\beta)$.  Since (\ref{datt}) agrees
with the Fuller Problem data for $(x,y)$ in some neighborhood
of $0$ and $|a|=1$, each $k> 0$ admits an $n(k)\in \N$ such that
$y^{\scr o}_{p(n(k))}(s,\gamma_{n(k)})=y^{\scr k}_{p(n(k))}
(s,\gamma_{n(k)})$ for all $s\ge 0$, so
\begin{eqnarray}
\int_0^\infty
\ell(y^k_{p(n(k))}(s,\gamma_{n(k)}),\gamma_{n(k)}(s))\, ds &=&
\sum_{n={n(k)}}^\infty\! \! \left[v_o(p(n))+\int_0^{1/(n+1)}
[s^2/2]^2\, ds\right]
 \nonumber\\ &\le&
\! \! \sum_{n={n(k)}}^\infty\! \! [M/n^4+1/(20\, n^5)] <
\infty,\nonumber
\end{eqnarray}
even though $y^k_{p(n(k))}(s,\gamma_{n(k)})\to (0,0)\not\in\T$ as
$s\to +\infty$.

One checks that $(H_1)$-$(H_6)$  hold for  (\ref{datt}) for {\em
all} $k\ge 0$.  For example, $(H_5)$ holds since the dynamics in
(\ref{datt}) agrees with the Fuller dynamics in a neighborhood of
the $y$-axis and  the Lagrangian $\ell$ assigns a positive cost to
staying at $(0,0)$ when $k>0$ and the Fuller Problem satisfies
$(H_5)$.     The fact that $STC\T$ holds for (\ref{datt}) follows
since $f(x,y,a)=(y-k,a)$ near $(k,k)$ and the Fuller Problem
satisfies $STC\{(0,0)\}$ (cf. \cite{M00}), along with a change of
coordinates.  Finally, condition $(H_6)$ holds by  Lemma
\ref{meagrelemma} with $g(x):=x^2$. This application of Lemma
\ref{meagrelemma} is based on the fact that $\Phi'$ has compact
support, which guarantees that the second derivative of the first
component of $y^k_x(s,\beta)$ is globally bounded.   We conclude
as follows: \bc{fullercor} Let $k\ge 0$ be constant, and choose
the exit time problem data (\ref{datt}).  If $w:\R^2\to\R$ is a
continuous function that is a bounded-from-below viscosity
solution of the corresponding HJBE
\[
\left[-y+k\Phi(x,y)\right] \left( Dw(x,y)\right)_1+|(Dw(x,y))_2|-x^2=0
\]
 on $\R^2\setminus\T$ that is
null at $\T$, then $w\equiv v_k$. \ec Taking $k=0$ in Corollary
\ref{fullercor} gives the uniqueness characterization for the
Fuller Problem HJBE asserted in \cite{MSCDC}. The novelty of
Corollary \ref{fullercor} is that it applies to problems violating
both the usual positivity condition (\ref{strong}) and  the
asymptotics condition (\ref{still}) from \cite{MSCDC}, and that it
establishes uniqueness of solutions of the HJBE in a class of
functions that includes functions that are not proper.

\br{newk} Using the fact that $x\mapsto x^2$ is convex, one shows
that $v_o$ is convex on $\R^2$ and therefore continuous. Moreover,
using Soravia's Backward Dynamic Programming Principle (cf.
\cite{BARD, Sor93a}), one can show that $(x,y)\mapsto
w(x,y):=-v_o(-x,y)$ is also a viscosity solution of the Fuller
Problem HJBE on $\R^2\setminus\{0\}$ vanishing at the origin.  The
argument is based on the facts that $v_o$ is a {\em bilateral}
viscosity solution of the HJBE and that each $p\in\R^2$ is an
optimal point (cf. \cite{BARD} for the definitions) and the fact
that $(p_1,p_2)\in D^+w(x,y)\implies (p_1,-p_2)\in D^-v_o(-x,y)$
and that $(p_1,p_2)\in D^-w(x,y)\implies (p_1,-p_2)\in
D^+v_o(-x,y)$. It follows that $v_o$ is the unique continuous
bounded-from-below viscosity solution of the corresponding HJBE on
$\R^2\setminus\{0\}$ that vanishes at the origin and that the
boundedness from below hypothesis of Corollary \ref{fullercor}
cannot be removed. \hfill\halmos\er

\br{mkk} Corollary \ref{fullercor} can be generalized.  For
example, the corollary remains true if the Lagrangian $\ell$ in
(\ref{datt}) is replaced by $\ell(x,y,a)=g(x)+k(1-|a|)^2$ for any
$g$ of class ${\cal MK}$, e.g., $g(x)=|x|^q$ for any $q>0$.  The
proof goes through without changes if the data are modified in
this way. Also, the target $\T=\{(k,k)\}$ can be replaced by
$\{(k,m)\}$ for any $k\ne 0$ and any $m\in\R$ if $\Phi$ is chosen
to be  $1$ near $(k,m)$ and zero in some open set containing the
$y$-axis.  Moreover, using the methods of $\S$\ref{do} below, the
above corollary can be extended to cover local and {\em
discontinuous} viscosity solutions. \hfill\halmos\er

\subsection{Degenerate Eikonal Equations}
\label{eikonalsubsection} This subsection shows how Theorem
\ref{thm1} applies to the HJBE's for a class of exit time problems
from geometric optics. The problems have the dynamics
$f(x,y,a,b)=(a,b)\in\overline{B_1(0)} \subseteq \R^2$ and the
Lagrangians
\begin{equation}
\label{wigglewiggle}
\ell(x,y,a,b)\; =\; \left[1+\sqrt{||(x,y)||}\right]^{-p},
\end{equation}
 where $p\ge 0$ is a constant that
will be further specify below.  (The argument we are about to give
also applies if we instead take the Lagrangian $(1+\sqrt{|x|})^p$
or $(1+\sqrt{|y|})^p$, or if the state space and compact control
set are in $\R^M$ for $M$ arbitrary.)

We choose any nonempty closed target $\T\subseteq\R^2$,
and we let $v_{e,p}$ denote the value function for the exit time
problem we have defined for each $p\ge 0$.  The corresponding HJBE
is
\begin{equation}
\label{eikeq}
||Dv(x,y)||=\frac{1}{\left[1+\sqrt{||(x,y)||}\right]^p},
\end{equation}
which is the eikonal equation of geometric optics for the propogation
of light in a medium with speed
\[c(x,y)=\left[1+\sqrt{||(x,y)||}\right]^p.\]
Viscosity solutions of eikonal equations have been
studied extensively (cf. \cite{BARD}, which covers cases where the
speed of the medium is bounded and also uniqueness questions for
eikonal equation solutions on bounded sets,
and \cite{S00}).  However, (\ref{eikeq})
is not covered by these results since $c$ is unbounded
and $\T$ may be unbounded.
It is  easy to check that for $0\le p\le 2$, the exit time problems
for these data satisfy $(H_1)$-$(H_6)$.  Indeed, if $q\in\R^2$
and if $\phi$ is any trajectory for $f$ starting at $q$, then we can
find a $K>0$ so that, for each $L>K$, we have
\[\ds
\int_0^L
\frac{ds}{\left[1+\sqrt{||\phi(s)||}\right]^p} \; \ge\;
\int_0^L
\frac{ds}{\left[1+\sqrt{||q||+s}\right]^p} \; \ge\;
\frac{1}{2}\int_K^L\frac{ds}{s^{p/2}}\; \to\; \infty
\]
as $L\to\infty$, so $(H_6)$ is satisfied vacuously.
We conclude as follows:

\bc{eikcor} Let $p\in[0,2]$ and  $\T\subseteq\R^2$ be closed and
nonempty.  If $w:\R^2\to\R$ is a continuous function that is a
bounded-from-below viscosity solution of (\ref{eikeq}) on
$\R^2\setminus\T$ that is null on $\T$, then $w\equiv v_{\scr
e,p}$. \ec

\br{rkkeik}
\label{ikk}
It was not necessary to assume that the target $\T$
is bounded.  If $p> 2$ in (\ref{wigglewiggle}), then Theorem \ref{thm1}
may not apply,
since $(H_6)$ could fail.  For example, if $p=4$,
and $\T=\{(x,0)\in\R^2: x\le
-1\}$  and $\beta\equiv (1,0)$, then (\ref{wigglewiggle}) gives
$\int_0^\infty \ell(y_{\scr
(0,0)}(s,\beta),\beta(s))\, ds<\infty$, even though   the trajectory
does not
remain bounded.
Moreover, the standard uniqueness
characterizations for exit time HJBE's (e.g., Corollary
IV.4.3 of \cite{BARD}) would not apply, since (\ref{strong}) is not
satisfied.  However, using
\cite{M01}, one can show that the statement of
Corollary \ref{eikcor} remains true {\em even without} the restriction
$p\in[0,2]$.  This is done by rewriting the HJBE (\ref{eikeq}) as
\begin{equation}
\label{neweik}
\left[1+\sqrt{||(x,y)||}\right]^p||Dv(x,y)||-1=0
\end{equation}
and then viewing (\ref{neweik}) as the HJBE for the exit time
problem with the non-Lipschitz dynamics
\[\tilde f(x,y,a,b)=\left[1+\sqrt{||(x,y)||}\right]^p(a,b)\]
 (with
$(a,b)\in\overline{B_1(0)}$
as before) and the Lagrangian $\tilde\ell\equiv 1$.  The dynamics
$\tilde f$ is then approximated by locally Lipschitz dynamics, and
then Theorem IV.4.4 of \cite{BARD} is applied.
For details,
see $\S$6.1
of \cite{M01}.\hfill\halmos
\er
\subsection{Shape-From-Shading Equations}
\label{shape}
Our results also apply to equations of the form
\[
I(x)\Psi(Du(x))-b(x)\cdot Du(x)-h^2(x)=0
\]
for $I$ nonnegative and $\Psi:\R^N\to\R$
any convex function with
$\Psi(0)=0$.  This
equation is studied in \cite{So99}.  Taking the Legendre transform
$\Psi^\star$ of $\Psi$, which is nonnegative, we can rewrite this equation
as
\[
\max_{a\in{\rm domain}(\Psi^\star)}
\left\{
-(b(x)-I(x)a)\cdot Du(x)-\left[h^2(x)+I(x)\Psi^\star(a)\right]\right\}=0.
\]
A particular case of this equation
(cf. \cite{So99})
is
\begin{equation}
\label{tone}
I(x)\left[1+||Du(x)||^2\right]^{1/2}-1=0,\; \; x\in\Omega\subseteq\R^2
\end{equation}
for open sets $\Omega$, which in fact
can be written as
\begin{equation}
\label{ttwo}
\ds
\max_{||a||\le 1}\left\{
I(x)a\cdot
Du(x)-\left[1-I(x)\left(1-||a||^2\right)^{1/2}\right]\right\}=0.
\end{equation}
The equation (\ref{ttwo}) arises in shape-from-shading models in
image processing, where $I(x)\in[0,1)$ is the intensity of light
reflected by an object (cf. \cite{S00}).  The objective in image
processing is to reconstruct the unknown function $u$,
representing the height of the surface on some subset $\Omega$ of
the plane, from the brightness of a single two-dimensional image
of the surface.  For the case of a Lambertian surface that is not
self-shadowing and that is illuminated by a single distant
vertical light source, the height $u$ is a viscosity solution of
(\ref{ttwo}).

Now pick any closed nonempty target
$\T\subseteq \R^2\setminus\{0\}$ and $\Omega:=\R^2\setminus\T$,
 and choose the intensity function
\begin{equation}
\label{pound0}
I(x):=\frac{||x||}{1+||x||}.
\end{equation}
Then (\ref{ttwo}) is an HJBE for an exit time problem with
the dynamics
\begin{equation}
\label{pound1}
f(x,u):=-I(x)u,\end{equation} the control set
$A=\overline{B_1(0)}\subseteq \R^2$,
and the Lagrangian
\begin{equation}
\label{pound2}
\ell(x,u)=1-\frac{||x||}{1+||x||}\left(1-||u||^2\right)^{1/2}.
\end{equation}
As explained in Remark \ref{rk2}, {}for general $\T\subseteq
\R^2\setminus\{0\}$, $\ell$ violates the positivity condition
(\ref{strong}) (since $\ell(x,0)\to 0$ as $||x||\to\infty$), so
the well-known results (e.g., those of \cite{BARD}) cannot be used
to get uniqueness characterizations for solutions of (\ref{tone}).
On the other hand, using the fact that \[||y_q(s,\beta)||\le
||q||+s\] for all $\beta\in {\cal A}$, $s\ge 0$, and $q\in \R^2$,
one can easily check that $(H_1)$-$(H_6)$ hold.  The argument is
similar to the validation of $(H_6)$ in
$\S$\ref{eikonalsubsection}. Therefore, we conclude from Theorem
\ref{thm1} that if $w:\R^2\to\R$ is a continuous function that is
a viscosity solution of (\ref{ttwo}) on $\R^2\setminus\T$ that
satisfies $(SC_w)$, then $w$ coincides with the shape-from-shading
value function. Local uniqueness characterizations and results for
discontinuous viscosity solutions  for the shape-from-shading
equation can also be given using the results in $\S$\ref{do}
below.

\br{nrk}
\label{nrkk}
As in the case of  eikonal equations, it was not necessary
to assume that the target $\T$ was bounded.  It is worth
remarking  that if we
replace the light intensity $I(x)$
with
\[
\tilde I(x):=\frac{3e^{2||x||}}{1+3e^{2||x||}}\; \; \in\; \; [3/4,1)
\]
in the previous example and keep the example the same otherwise, then
Theorem \ref{thm1}  would no longer apply, since condition $(H_6)$ may not
be
satisfied.  However, for such cases, we can still apply
\cite{M00} to get uniqueness of {\em proper} solutions of the
corresponding HJBE.
For example, take ${\cal T}=\{(0,r): r\le -1\}$ and the
control
\[
\beta(t)\equiv (0,-1/(t+1)),\]
and let $t \mapsto x(t)=(x_1(t),x_2(t))$ denote
the trajectory
of
$
\tilde f(p,u):=-\tilde I(p)u$ for the initial position $p(0)=(0,1)$ and
the control $u=\beta(t).
$
For all $t>0$, we then have $x_1(t)=0$,
\[
\ds x_2(t)\; \; =\; \;  1+\int_0^t\tilde I(x(s))\frac{1}{s+1}\, ds
\; \;
\ge\; \;
1+\frac{3}{4}\int_0^t\frac{1}{s+1}\, ds\; \; =\; \;
1+\frac{3}{4}\ln(t+1),\]
so $x(t)$ is not bounded.  However,
\begin{eqnarray}
\label{calc} e^{2||x(t)||}||\beta(t)||^2& = &
\frac{1}{(t+1)^2}e^{\left[2+2\int_0^t\tilde I(x(s))\frac{1}{s+1}\,
ds\right]}\nonumber \\&\le&
\frac{1}{(t+1)^2}e^2e^{2\ln(t+1)}\nonumber\\&\le&  e^2.
\end{eqnarray}

Therefore,  if $\tilde \ell(p,u)=1-\tilde
I(p)[1-||u||^2]^{1/2}$
denotes the corresponding Lagrangian, then since
we have \[\tilde \ell(p,u)\le 1-\tilde I(p)[1-||u||^2]\; \; \forall p\in
\R^2,\; u\in \overline{B_1(0)},\] (\ref{calc}) gives
\begin{eqnarray}
\ds\int_0^\infty \tilde \ell(x(s),\beta(s))\, ds &\le &
\int_0^\infty \frac{1+3e^{2||x(s)||}||\beta(s)||^2}{1+3e^{2||x(s)||}}\, ds
\nonumber\\
&\le &
\left[1+3e^2\right]
\int_0^\infty\frac{dt}{1+3e^{2+3/2\ln(t+1)}}\; \; \; \; \;
\nonumber\\
&\le & \left(\left[1+3e^2\right]/\left[3e^2\right]\right)
\int_0^\infty\frac{dt}{(t+1)^{3/2}}\; <\; \infty,\nonumber
\end{eqnarray}
even though $t\mapsto x(t)$
is not bounded, which shows $(H_6)$ is not satisfied.
Moreover, the standard uniqueness characterizations for exit time HJBE's
(cf. \cite{BARD, BS92}) would again not apply,
since the Lagrangian
$\tilde \ell$ is not uniformly bounded below by positive constants.
However, since $(H_5)$ holds,
 one can use
\cite{M99a} to show that for any nonempty closed
target $\T\subseteq \R^2\setminus\{0\}$, any
proper continuous viscosity solution
of the
corresponding HJBE
\[
\ds \sup_{||a||\le 1}
\left\{
\tilde I(x)a\cdot Du(x)-\left[1-\tilde
I(x)\left(1-||a||^2\right)^{1/2}\right]\right\}\; \; =\; \; 0
\]
on $\R^2\setminus {\cal T}$ that is null on ${\cal T}$ must in
fact be identically equal to the shape-from-shading exit time
value function $v_{\scr sfs}$ for the target ${\cal T}$, the
dynamics $\tilde f$, and the Lagrangian $\tilde
\ell$.\hfill\halmos \er

\br{at} Notice that it was not necessary to assume that the domain
set $\Omega$  for (\ref{tone}) was bounded. It is worth pointing
out that one cannot in general expect uniqueness of solutions for
the shape-from-shading HJBE for cases where $I$ is allowed to take
the value $1$, since  the surface $u$ and $-u$ could both be
viscosity solutions of (\ref{tone}).  For example, take the light
intensity $I(x)=(1+4||x||^2)^{-1/2}$, $\T=\R^2\setminus B_1(0)$,
and the surface $u(x)=1-||x||^2$ on $B_1(0)$ and zero elsewhere.
Clearly, $u$ and $-u$ are both solutions of (\ref{tone}). However,
$(H_5)$-$(H_6)$ are not satisfied, since the trajectory
$\phi(t)\equiv 0$ gives zero integrated costs on $[0,\infty)$
without ever reaching the target, so this case is not covered by
Theorem \ref{thm1}. For the analysis of cases where
$\#\{x:I(x)=1\}=1$, see \cite{LRT}, and for bounded viscosity
solutions of (\ref{ttwo}), see \cite{RT}.\hfill\halmos \er

\section{Discontinuous and Local HJBE Solutions}
\label{do}
This section gives variants of Theorem \ref{thm1} for discontinuous
and local HJBE
solutions.
We study
discontinuous solutions using the
envelopes approach from \cite{BARD}.

\subsection{A Remark on  Discontinuous Viscosity Solutions}
\label{disco} Under $(H_1)$-$(H_6)$, the value function $v_g$
could be  discontinuous (cf. \cite{BARD}, pp. 248-249). This
suggests the question of how one can characterize $v$ as the
unique {\em discontinuous}  solution of the HJBE on
$\R^N\setminus\T$ that satisfies $(SC_v)$. By a discontinuous
solution, we mean the following.  For each locally bounded
function $w: S\to\R$ on a set $S\subseteq\R^N$, we define the
following semicontinuous envelopes:
\[
\ds
w_\star(x):=\liminf_{S\ni y\to x} w(y) \; \; \text{and}\; \;
w^\star(x):=\limsup_{S\ni y\to x} w(y)\; \; .
\]
We call $w_\star$ the {\bf lower envelope}  of $w$, and we call
$w^\star$ the {\bf upper envelope}  of $w$.  For $\cg$, $S$, and
$F$ satisfying the requirements of Definition \ref{definitionvis},
we then say that a locally bounded  function $w:S\to\R$ is a {\bf
discontinuous subsolution} (resp., {\bf supersolution}) of
$F(x,Dw(x))=0$ on $\cg$ provided $F(x_o, D\gamma(x_o))\le 0$
(resp., $\ge 0$) for each $\gamma\in C^1(\cg)$ and each local
maximizer (resp., minimizer) of $w^\star-\gamma$ (resp.,
$w_\star-\gamma$) on $\cg$.\footnote{In this context,
`discontinuous' means ``not necessarily continuous''.} A {\bf
(discontinuous viscosity) solution} of $F(x,Dw(x))=0$ on $\cg$ is
then a function that is simultaneously a discontinuous subsolution
and a discontinuous supersolution of $F(x,Dw(x))=0$ on $\cg$.
Lemma \ref{lemma0} remains true if $u\in C(\bar E)$ is replaced by
any bounded discontinuous subsolution of the HJBE on $E$ and $u$
in (\ref{halff}) is replaced by $u^\star$. Also, Lemma
\ref{lemma1} remains true if $w\in C(\bar B)$ is replaced by any
bounded discontinuous supersolution of the HJBE on $B$ and $w$ in
(\ref{quadruple}) is replaced by $w_\star$.  Using these facts,
one can prove the following generalization of Theorem \ref{thm1}:
If $(H_1)$-$(H_6)$ hold, if $w:\R^N\to\R$ is a discontinuous
viscosity solution of the HJBE on $\R^N\setminus\T$ that satisfies
$(SC_w)$, and if $w_\star$ is continuous on $\R^N$, then $w\equiv
v$ on $\R^N$. \footnote{The continuity of $w_\star$ is used to
ensure that the sets $S_\kappa$ in the proof of Theorem \ref{thm1}
are open.  The condition that $w_\star$ is continuous of course
holds automatically if $w$ is continuous. However, $(SC_w)$ and
continuity of $w_\star$ can even be satisfied by functions that
are nowhere continuous.  For example, if we take the indicator
function $w\equiv \ind_{{\mathbb Q}}:\R\to\{0,1\}$, then
$w_\star\equiv 0$.  This generalized version of Theorem \ref{thm1}
remains true if the pointwise condition that $w\equiv 0$ on $\T$
is replaced by the less restrictive requirement that there be a
locally bounded function $g:\R^N\to\R$ for that
\[
\forall x\in\T, \; \; w_\star(x)\ge g_\star(x) \; \text{and}\;
w^\star(x)\le g^\star(x)
\]
except that the conclusion that $w\equiv v$ is replaced by the
following inequalities on $\R^N$ (cf. Remark \ref{rk1}):
$w_\star\ge v_{g_\star}$ and $w^\star\le v_{g^\star}$. In case
$g\in C(\T)$, this implies $w\equiv v_g$ on $\R^N$.} The proof is
almost identical to the proof of Theorem \ref{thm1} but with $w$
replaced by $w_\star$ in the proof of the inequality $w\ge v$. For
cases where $v$ is continuous, this establishes that all solutions
$w$  of the HJBE (\ref{hjb}) on $\R^N\setminus\T$ that satisfy
$(SC_w)$ and continuity of $w_\star$ agree with $v$, and therefore
are continuous.

\subsection{Local Solutions of the HJBE}
\label{obsies}
This subsection shows how to extend Theorem \ref{thm1} to get
uniqueness of solutions of the HJBE on sets of the form
$\Omega\setminus\T$ for open sets $\Omega$.  We set
\[ \crr=\left\{x\in\R^N: \inf\{t_x(\beta): \beta\in
{\cal A}\}<\infty\right\}, \]
so $\crr$ is the set of points that can be brought
to $\T$ in finite time using the dynamics $f$.  Using $(H_2)$-$(H_3)$, one
shows that $\crr$ is open (cf. \cite{BARD}).  In  many classical cases
where
$\ell$ is
bounded below by a positive constant, one has \begin{equation}
\label{vco} v \text{is\ continuous\ on} \crr, \; \;
\text{and}\; \;
\lim_{x\to x_o}v(x)=+\infty\;  \forall x_o\in \partial \crr.
\end{equation} On the other hand, one easily finds examples where
$\ell$ is not bounded below by a positive constant and the limit
condition in (\ref{vco}) fails.  Here is an elementary example
where this occurs: \begin{example}\label{bigg} Take $N=1$,
$\T=[1,+\infty)$, $A=\{+1\}$, $f(x,a)=|x|a$, and $\ell(x,a)=|x|$.
In this case, \[ v(\bar x)= \int_0^{\ln(1/\bar x)}\bar xe^t\,
dt=1-\bar x\to 1\; \; \text{as}\; \; \bar x\downarrow 0,\] even
though $0\in\partial \crr$. \hfill\halmos \eex

This motivates the question of how one can characterize $v$ as a
unique viscosity solution of the HJBE on $\crr\setminus\T$ for
cases where $\crr\ne \R^N$ and the extra condition (\ref{vco})
holds.  To address this question, we assume the following relaxed
version \bi\item[]\bi \item[$(H'_5)$]\ \ If $x\in\crr\setminus\T$
and $\beta\in {\cal A}$, then
$\int_0^t\ell^r(y_x(s,\beta),\beta(s))\, ds>0$ for all $t\in(0,
\infty)$. \ei\ei of $(H_5)$. We also fix an open set
$\Omega\subseteq \crr$ containing $\T$, and we consider viscosity
solutions of the HJBE (\ref{hjb}) on $\Omega\setminus\T$ that
satisfy the localization \bi\item[]\bi\item[]\bi \item[$(OSC_{\scr
w,\Omega})$] $w$ is bounded-from-below on $\Omega$, $w\equiv 0$ on
$\T$, and $\ds\lim_{x\to x_o}w(x)=+\infty$ $\forall x_o\in
\partial\Omega$.\ei\ei\ei Noting that $v$ satisfies $(OSC_{\scr
v,\crr})$ if (\ref{vco}) holds, we then have the following local
version of Theorem \ref{thm1}: \bt{thm5} Let $(H_1)$-$(H_4)$,
$(H'_5)$, and $(H_6)$ hold. Let $\Omega\subset\crr$ be an open set
containing $\T$.  Let $w:\Omega\to\R$ be a continuous function
that is viscosity solution of the HJBE (\ref{hjb}) on
$\Omega\setminus\T$ that satisfies $(OSC_{\scr w, \Omega})$. Then,
$w\equiv v$ on $\Omega$. In particular, if $v$ satisfies
(\ref{vco}), then $v$ is the unique viscosity solution $w$ of the
HJBE on $\crr\setminus\T$ in the class of all continuous functions
$w:\crr\to\R$ that satisfy $(OSC_{\scr w, \crr})$.\et

\br{prff} The proof of the inequality $w\le v$ for Theorem
\ref{thm5} is exactly the proof of that inequality in \cite{M00}.
The proof is slightly more complicated than the proof that $w\le
v$ for Theorem \ref{thm1}, since one must consider trajectories
that reach $\T$ in finite time but that exit $\Omega$ before the
first  time they ever reach $\T$. The proof of the reverse
inequality closely follows the proof of Theorem \ref{thm1} except
that instead of setting ${\cal S}={\cal S}_\kappa\cap B_J(0)$, we
set ${\cal S} = {\cal S}_\kappa\cap B_J(0)\cap \Omega$. We rule
out cases where $\bar x_J\in \partial\Omega$ using the limit
condition in $(OSC_{\scr w,\Omega})$.   Theorem \ref{thm5} can
also be generalized to the case of discontinuous viscosity
solutions using the method of $\S$\ref{disco}.\er

\section{Problems with Unbounded Control Sets}
\label{noncompact} We close by giving two variants of Theorem
\ref{thm1} that can be applied for cases where the control set
$A\subseteq \R^M$ is closed but possibly unbounded. In the first
variant, we impose regularity conditions on the data that penalize
the use of control set values of large norm.  In the second
variant, we replace the possibly unbounded control set $A$ with a
suitable compact set of vector field valued controls. Recall the
definition (\ref{noncom}) of ${\cal A}$ that applies to possibly
noncompact control sets.
\subsection{Penalization Method}
For simplicity, let us assume that all the sets \[{\cal
D}(x):=\{(f(x,a),\ell(x,a)): a\in A\}\] are convex.
As explained in $\S$\ref{dh},
the set of inputs
$\alpha\in {\cal A}$  can then be taken to be
the measurable functions valued in $A$ (by the Filippov Selection
Theorem).  We
assume that
$(H_2)$-$(H_6)$ are satisfied, where
$0\in A\subseteq\R^M$ for $M\in \N$
and $A$ is closed
but not necessarily compact.  Following \cite{BDL97, DAL, M00},
we then add the following conditions on $f$ and $\ell$:

\bi\item[]\bi \item[$(H_7)$] $f$ is bounded on $B_R(0)\times A$
for each $R>0$. \item[$(H_8)$] There is a modulus $\omega$ such
that $|\ell(x,u)-\ell(y,u)|\le \omega(||x-y||)$ for all
$x,y\in\R^N$ and $u \in A$. \item[$(H_9)$] There exist constants
$\ell_o>0$, $C_o\ge 0$, $\beta\in(0,1]$, $\delta_2\ge 0$,
$\bar\ell\ge 0$, and $\delta_1> 1$ such that the following
conditions hold for all $x,y\in\R^N$ and $a\in A$: \bi \item[(a)]
$\ell(x,a)\ge \ell_o||a||^{\delta_1}-C_o$ \item[(b)]
$|\ell(x,a)-\ell(y,a)|\le
\bar\ell||x-y||^\beta(1+||a||^{\delta_1}+
||x||^{\delta_2}+||y||^{\delta_2})$ \ei\ei\ei (Recall that a {\bf
modulus} is a nondecreasing continuous function
$\omega:[0,\infty)\to[0,\infty)$ for which $\omega(0)=0$.) As
shown in \cite{BARD}, Lemmas \ref{lemma0} and \ref{lemma1} remain
true if $(H_2)$-$(H_8)$ are assumed instead of the assumptions
$(H_1)$-$(H_6)$. These assumptions penalize the use of control set
values of large norm. We then consider only viscosity solutions
$w$ of the HJBE on $\R^N\setminus\T$ for which the subdifferential
sets $D^-w(x)$ are locally bounded, i.e., such that $\sup\{||p||:
p\in  D^-w(x), x\in K\}<\infty$ for each compact set
$K\subseteq\R^N$. As shown in Theorem I.7.3 of \cite{CLSW98}, this
is equivalent to considering only locally Lipschitz solutions of
the HJBE on $\R^N\setminus\T$.
  In this case,  the infimizations in the
restriction of the HJBE
to any $B_J(0)$
can be taken over a corresponding compact set
$C_J\subset A$, i.e., in the notation we introduced in $\S$\ref{dh},
$H_{A}\lceil [B_J(0)\times D_J] =
H_{C_J}\lceil [B_J(0)\times D_J]$, where $D_J$ is a bounded set large
enough to contain
$\{p\in D^-w(x): x\in B_J(0)\}$
(cf.
\cite{BDL97, DAL} for the proof).
Then the
arguments in $\S$\ref{prf} on $B_J(0)$ apply with the compact control set
$C_J$
replacing
$A$, and then we iterate on $J$ to
get an input $\hat\alpha:[0,\infty)\to A$ as before.    We
then invoke $(H_6)$ to conclude as
follows:

\bt{thm2} Assume hypotheses $(H_2)$-$(H_9)$, with $A$ a closed set
containing $0\in \R^M$.    Let $w:\R^N\to \R$ be a locally
Lipschitz  function that is a viscosity solution of (\ref{hjb}) on
$\R^N\setminus\T$ that satisfies $(SC_w)$.  Then $w\equiv v$. \et

\subsection{Vector Field Valued Controls Method}
Another way to extend Theorem \ref{thm1} to the case of noncompact control
sets is as follows.  As in the previous subsection, we assume the sets
${\cal D}(x)$ are all
convex.  We give $C(\R^N,
\R^N\times \R)$ the topology of compact convergence (cf. \cite{Munkres}).
We continue to assume $(H_2)$-$(H_8)$ and that $A\subseteq \R^M$
is
closed and nonempty but possibly unbounded.  We also add the
following assumptions:
\bi\item[]\bi
\item[$(NC_1)$]
\ \ $\sup\{\ell(0,u):u\in A\}\; <\; \infty$.
\item[$(NC_2)$]
\ \ $\{(f(\cdot, u),\ell(\cdot, u)):u\in A\}\; \subseteq \;
C(\R^N, \R^N\times \R)$ is closed.
\ei\ei
These guarantee that the supremum in the definition of the HJBE is always
finite.
It
follows from the Ascoli-Arzel\'a Theorem that
$
K := \{k_u(\cdot) := (f(\cdot,u), \ell(\cdot, u)): u\in
A\}
$
is a compact subset of the metric space $C(\R^N, \R^N\times \R)$
(cf. \cite{M00,Munkres}).
Define the projection mappings $\pi_j$ on $K$ by
\[
\pi_j(k_u(\cdot))=\twoif{f(\cdot, u),}{j=1}{ \ell(\cdot,
u),}{j=2}\; \; \; \forall u\in A.
\]
We now apply the method of our proofs to the new exit time problem
whose dynamics $F$, Lagrangian $\Lambda$, and set $\tilde{\cal A}$
of admissible controls
 are
\[
F(x,k) = (\pi_1\circ k)(x), \; \; \Lambda(x,k)= (\pi_2\circ
k)(x)\;  \& \; \tilde{\cal A} := \{[0,\infty)\ni t\mapsto
k_{\beta(t)}: \beta\in {\cal A}\}\] with the same target $\T$.
Notice that $\{(F(x,k), \Lambda(x,k)): k\in K\}$ is convex for
each $x\in \R^N$. Let $\tilde v$ denote the value function of this
new problem. Since the trajectories of $F$ with the controls
$\tilde {\cal A}$ are exactly the trajectories of $f$ with
controls in $\cal A$, it follows that $\tilde v\equiv v$.
Moreover, the new problem satisfies $(H_1)$-$(H_6)$ (with $K$
replacing $A$, $F$ replacing $f$, and $\Lambda$ replacing $\ell$).
Our proof of Theorem \ref{thm1} then gives the following:
\bt{thm3} Let $\emptyset\ne A\subseteq \R^M$ be closed.  Assume
$(H_2)$-$(H_8)$ and $(NC_1)$-$(NC_2)$. Let $w:\R^N\to \R$ be a
continuous function that is a viscosity solution of (\ref{hjb}) on
$\R^N\setminus\T$ that satisfies $(SC_w)$.  Then $w\equiv v$. \et
We remark that if $(H_2)$-$(H_8)$ and $(NC_1)$-$(NC_2)$ all hold
with $A\ne\emptyset$ a closed subset of $\R^N$, and if
$\crr=\R^N$, then the value function $v$ is a discontinuous
viscosity solution of the HJBE on $\R^N\setminus\T$ (cf.
\cite{BARD}).  If we also assume $v_\star$ is continuous, then a
generalization of Theorem \ref{thm3} characterizes $v$ as the
unique discontinuous viscosity solution $w$ of the HJBE in the
class of functions $w:\R^N\to\R$ that satisfy $(SC_w)$ and
continuity of $w_\star$. The generalization of Theorem \ref{thm3}
to discontinuous solutions follows from the argument of
$\S$\ref{disco}.  Also, the theorem extends to local HJBE
solutions using the arguments of the previous section.

\end{document}